\documentclass{ifacconf}

\usepackage{graphicx}      
\usepackage{natbib}        
\usepackage{amssymb}
\usepackage{amsmath}

\begin{document}
\begin{frontmatter}

\title{More Details on Analysis of Fractional-Order Lotka-Volterra Equation}


\author[First]{F. Merrikh-Bayat}

\address[First]{Department of Electrical and Computer Engineering, University of Zanjan, Zanjan, IRAN (e-mail: f.bayat@znu.ac.ir)}

\begin{abstract}                
According to the long-memory principle appears in fractional-order
dynamical systems, analysis of these systems is commonly more
complicated than those described by nonlinear ordinary
differential equations. Another difficulty is due to the fact that
some classical tools such as the Lyapunov stability theorems
cannot directly be applied to nonlinear fractional differential
systems. The aim of this paper is to study the stability of a
general form of the two-dimensional fractional-order
Lotka-Volterra equation and numerically investigate the domain of
attraction of its stable equilibrium points. It is also shown that
the fractional-order Lotka-Volterra system can never have a stable
focus although the linearized system has complex conjugate modes.
Some other properties of the fractional-order Lotka-Volterra
equation, which are not observed in integer-order case, are also
discussed.
\end{abstract}


\end{frontmatter}
\section{Introduction}
The simplest model of predator-prey interactions was first
developed independently by  Alfred J. Lotka \citep{lotka} and Vito
Volterra \citep{volterra}. The classical two-dimensional
Lotka-Volterra equation is given by:
\begin{equation}
\left\{%
\begin{array}{ll}
    \frac{dy_1}{dt}=ay_1-by_1y_2, \\
    \frac{dy_2}{dt}=-cy_2+by_1y_2, \\
\end{array}%
\right.
\end{equation}
where $y_2$ is the number of some predator (for example, wolves),
$y_1$ is the number of its prey (for example, rabbits) and $a$,
$b$, and $c$ are real parameters representing the interaction of
the two species. This system and its extensions have been fully
studied before by several researchers \citep{murray}.

Recently, \cite{ahmed07} introduced the fractional-order
Lotka-Volterra predator-prey system:
\begin{equation}\label{lotka_ahmed}
\left\{%
\begin{array}{ll}
    D^{\alpha}x_1(t)=x_1(t)(r-ax_1(t)-bx_2(t)), \\
    D^{\alpha}x_2(t)=x_2(t)(-d+cx_1(t)), \\
\end{array}%
\right.
\end{equation}
for some positive real constants $a,~b,~c,~d$, and
$\alpha\in(0,1]$, and studied the stability of its equilibrium
points (see Section \ref{sec_back} for the definition of
fractional-order operators).

The aim of this paper is to study a more general form of the
fractional-order Lotka-Volterra equation (as defined in
(\ref{lotka1})) in which the order of fractional derivatives are
assumed to be different in general. Such a model may be used for
better modelling and understanding the behavior of more
complicated predator-pray systems.

Currently, nonlinear fractional-order systems constitute a
challenging research area, mainly because some powerful classical
tools such as the Lyapunov's stability method cannot directly be
applied to these systems. It is also a well-understood fact that
determining the domain of attraction of stable equilibrium points
of nonlinear fractional-order systems is not a straightforward
task. The main reason for this difficulty is that the phase
portrait of these systems cannot be plotted using the so-called
isoclines. The long-memory principle of the fractional-order
operators leads to some strange behaviors which are also discussed
in this paper.

The rest of this paper is organized as follows. Some mathematical
preliminaries are reviewed in Section \ref{sec_back}. Stability of
the fractional-order Lotka-Volterra equation under consideration
is studied in Section \ref{sec_stab}. Section \ref{sec_domain}
discusses on the domain of attraction of the stable equilibrium
points of the fractional-order Lotka-Volterra equation introduced
in Section \ref{sec_stab}. Some properties of the fractional-order
Lotka-Volterra equation are studied in Section 5, and finally
Section 6 concludes the paper.

\section{Mathematical background}\label{sec_back}
\subsection{fractional-order operators}
Three kinds of definitions are widely used to define
fractional-order derivatives: Gr\"{u}nwald-Letnikov derivative,
Riemann-Liouville derivative and Caputo derivative
\citep{podlubny99}. These three definitions are in general not
equivalent. The definition of the fractional derivative given by
Caputo has the advantage of only requiring initial conditions
given in terms of integer-order derivatives. Clearly, such initial
conditions represent well-understood features of a physical
situation. In this paper, we will also use the Caputo derivative
because of its applicability to real world models.

The Caputo fractional derivative of order $\alpha$ of function
$f$, which is denoted by $_a^CD_t^\alpha f(t)$, is defined as
\begin{equation}\small{
_a^CD_t^\alpha f(t)=\left\{%
\begin{array}{ll}
    \frac{1}{\Gamma(n-\alpha)}\int_a^t\frac{f^{(n)}(x)}{(t-x)^{\alpha-n+1}}dx,~~~n-1<\alpha<n\\
    f^{(n)}(t),\hspace{4cm} \alpha=n\\
\end{array},%
\right.}
\end{equation}
where $t>a$ and $n\in\mathbb{Z}^+$. For simplicity of the
notation, the symbol $D^\alpha$ is used in the rest of this paper
to indicate $_0^CD_t^\alpha$.
\subsection{stability of fractional-order
systems}\label{subsec_stab} The following theorem will be
instrumental in what follows.

\begin{thm} \citep{deng07} Consider the $n$-dimensional
linear fractional-order system:
\begin{equation}\label{system}
\left\{%
\begin{array}{ll}
    D^{\alpha_1}x_1=a_{11}x_1+a_{12}x_2+\ldots+a_{1n}x_n\\
    D^{\alpha_2}x_1=a_{21}x_1+a_{22}x_2+\ldots+a_{2n}x_n\\
    \vdots\\
    D^{\alpha_n}x_n=a_{n1}x_1+a_{n2}x_2+\ldots+a_{nn}x_n\\
\end{array},%
\right.
\end{equation}
where all $\alpha_i$'s are rational numbers between $0$ and $1$.
Assume $M$ be the lowest common multiple of the denominators
$u_i$'s of $\alpha_i$'s, where $\alpha_i=v_i/u_i$, $(u_i,v_i)=1$,
$u_i,v_i\in\mathbb{Z}^+$, for $i=1,2,\ldots,n$. Define
\begin{equation}
\Delta(\lambda)=\left(%
\begin{array}{cccc}
  \lambda^{M\alpha_1}-a_{11} & -a_{12} & \ldots & -a_{1n} \\
  -a_{21} & \lambda^{M\alpha_2}-a_{22} & \ldots & -a_{2n} \\
  \vdots & \vdots & \ddots & \vdots \\
  -a_{n1} & -a_{n2} & \ldots & \lambda^{M\alpha_n}-a_{nn} \\
\end{array}%
\right).
\end{equation}
Then system (\ref{system}) is globally asymptotically stable (in
the Lyapunov sense) if all roots $\lambda$'s of the equation $\det
(\Delta (\lambda))=0$ satisfy $| \arg(\lambda)|>\pi/(2M)$.
\end{thm}


It is not difficult to show that the condition given in Theorem 1
is the necessary and sufficient condition for asymptotic stability
of system (\ref{system}) (see also \citep{matignon} for more
details on this subject). Now, consider the nonlinear
fractional-order system:
\begin{equation}\label{system2}
D^{\alpha_i}x_i(t)=f_i(x_1,x_2,\ldots,x_n),\quad i=1,2,\ldots,n,
\end{equation}
where all $\alpha_i$'s are rational numbers between 0 and 1.
Clearly, equilibrium points of (\ref{system2}) are roots of the
equation
\begin{equation}
f_i(x_1,x_2,\ldots,x_n)=0,\quad i=1,2,\ldots,n.
\end{equation}
Now, let $\mathbf{x}^*=(x_1^*,x_2^*,\ldots,x_n^*)$ be an
equilibrium point of (\ref{system2}), i.e.
$f_i(x_1^*,x_2^*,\ldots,x_n^*)=0$ for $i=1,2,\ldots,n$. Suppose
that $M$ is the lowest common multiple of the denominators $u_i$'s
of $\alpha_i$'s, where $\alpha_i=v_i/u_i$, $(u_i,v_i)=1$,
$u_i,v_i\in\mathbb{Z}^+$, for $i=1,2,\ldots,n$. Then, according to
Theorem 1 it can be easily shown that $\mathbf{x}^*$ is
asymptotically stable if and only if the inequality:
\begin{equation}
|\arg(\lambda)|>\pi/2M,
\end{equation}
holds for all roots $\lambda$'s of the equation:
\begin{equation}\label{det_eq}
\det \left({\rm diag}([\lambda^{M\alpha_1}, \lambda^{M\alpha_2},
\ldots, \lambda^{M\alpha_n}])-\mathbf{J} \right)=0,
\end{equation}
where the notation ${\rm diag}([r_1, r_2, \ldots, r_n])$ denotes
an $n\times n$ diagonal matrix as follows:
\begin{equation}
{\rm diag}([r_1, r_2, \ldots, r_n])=\left(%
\begin{array}{cccc}
  r_1 & 0 & \ldots & 0 \\
  0 & r_2 & \ldots & 0 \\
  \vdots & \vdots & \ddots & \vdots \\
  0 & 0 & \ldots & r_n \\
\end{array}%
\right),
\end{equation}
and
\begin{equation}
\mathbf{J}=\partial \mathbf{f}/\partial \mathbf{x}|_{x^*},\quad
\mathbf{f}=[f_1\quad f_2\quad \ldots \quad f_n]^T.
\end{equation}

\section{Stability analysis of the two-dimensional fractional-order Lotka-Volterra
equation}\label{sec_stab} In this paper, we mainly study the
two-dimensional fractional-order Lotka-Volterra equation:
\begin{equation}\label{lotka1}
\left\{%
\begin{array}{ll}
    D^\alpha y_1=f_1(y_1,y_2)=y_1(a-by_2), \\
    D^\beta y_2=f_2(y_1,y_2)=y_2(-c+by_1), \\
\end{array}%
\right.
\end{equation}
where $\alpha$ and $\beta$ are positive rational constants. Note
that it is not a considerable loss of generality to limit the
studies to the case that both $\alpha$ and $\beta$ are rational
numbers since, in practice, all numbers are stored with a limited
precision in computer and moreover, one can find a rational number
in any neighborhood of a given nonrational number. In the
following, we discuss on the asymptotic stability of the
equilibrium points of (\ref{lotka1}) in two cases separately.

\subsection{Case 1: $\alpha$ and $\beta$ between 0 and 1}
The stability theorem presented in Section \ref{subsec_stab} can
directly be used to study the stability of the equilibrium points
of (\ref{lotka1}) when both $\alpha$ and $\beta$ are rational
numbers in the range $(0,1)$. In this case the system has two
equilibrium points denoted as $\mathbf{Y}_1^*=(0,0)$ and
$\mathbf{Y}_2^*=(c/b,a/b)$. The Jacobian matrices are calculated
as
\begin{equation}
\mathbf{J}_1=\left(%
\begin{array}{cc}
  a & 0 \\
  0 & -c \\
\end{array}%
\right),
\end{equation}
and
\begin{equation}
\mathbf{J}_2=\left(%
\begin{array}{cc}
  0 & -c \\
  a & 0 \\
\end{array}%
\right),
\end{equation}
respectively at $\mathbf{Y}_1^*$ and $\mathbf{Y}_2^*$. Assume that
$\alpha=u/M$ and $\beta=v/M$ for some positive integers $u$, $v$,
and $M$. At $\mathbf{Y}_1^*$, (\ref{det_eq}) concludes that
\begin{equation}\label{charac1}
\det \left( \left(%
\begin{array}{cc}
  \lambda^u & 0 \\
  0 & \lambda^v \\
\end{array}\right)%
-
\left(%
\begin{array}{cc}
  a & 0 \\
  0 & -c \\
\end{array}%
\right) \right)=\left(\lambda^u-a\right)
\left(\lambda^v+c\right)=0.
\end{equation}
So, the equilibrium point $\mathbf{Y}_1^*$ is stable if and only
if all roots of (\ref{charac1}) lie in the sector of stability
defined by
\begin{equation}\label{sector_cond}
|\arg(\lambda)|>\frac{\pi}{2M}.
\end{equation}
But if $a>0$ or $c<0$ then at least one of the $u+v$ roots of
(\ref{charac1}) lie in the sector defined by
$|\arg(\lambda)|<\pi/(2M)$ and consequently, the equilibrium point
$\mathbf{Y}_1^*$ becomes unstable. Hence, the necessary condition
for stability of $\mathbf{Y}_1^*$ is to have $a<0$ and $c>0$.
Assuming $a<0$ and $c>0$, the roots of (\ref{charac1}) are
calculated as
\begin{equation}
\lambda=\sqrt[u]{|a|}e^{i(2h+1)\pi/u},\quad h=0,1,\ldots,u-1,
\end{equation}
and
\begin{equation}
\lambda=\sqrt[v]{c}e^{i(2h+1)\pi/v},\quad h=0,1,\ldots,v-1,
\end{equation}
where $i=\sqrt{-1}$. Clearly, all of the above roots lie in the
sector of stability defined by (\ref{sector_cond}) if and only if
we have
\begin{equation}
\frac{\pi}{u}>\frac{\pi}{2M}\quad \mathrm{and}\quad
\frac{\pi}{v}>\frac{\pi}{2M},
\end{equation}
which concludes that
\begin{equation}
\alpha<2\quad \mathrm{and}\quad \beta<2,
\end{equation}
considering the fact that $\alpha=u/M$ and $\beta=v/M$. As a
result, $\mathbf{Y}_1^*$ is a stable equilibrium point of
(\ref{lotka1}) if and only if we have $a<0$ and $c>0$ provided
that $\alpha$ and $\beta$ are real numbers between 0 and 1. Note
that the stability of $\mathbf{Y}_1^*$ is independent of the value
assigned to $b$.

\begin{figure}
\begin{center}
\includegraphics{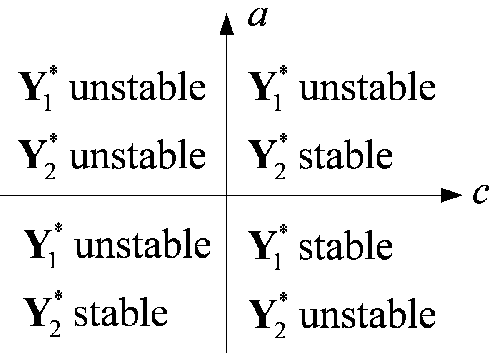}
\caption{Stability properties of the equilibrium points of
(\ref{lotka1}) in $a-c$ plane when both $\alpha$ and $\beta$ are
rational numbers between 0 and 1.}\label{fig_reg}
\end{center}
\end{figure}

At $\mathbf{Y}_2^*$, (\ref{det_eq}) reads
\begin{equation}\label{charac2}
\det \left( \left(%
\begin{array}{cc}
  \lambda^u & 0 \\
  0 & \lambda^v \\
\end{array}\right)%
-
\left(%
\begin{array}{cc}
  0 & -c \\
  a & 0 \\
\end{array}%
\right) \right)=\lambda^{u+v}+ac=0.
\end{equation}
Clearly, if $a$ and $c$ satisfy the inequality $ac<0$ then
(\ref{charac2}) will have at least one root outside the region of
stability defined by (\ref{sector_cond}) and hence,
$\mathbf{Y}_2^*$ will be unstable. So, the necessary condition for
the stability of $\mathbf{Y}_2^*$ is that we have $ac>0$. Assuming
$ac>0$, (\ref{charac2}) yields
\begin{equation}
\lambda=\sqrt[u+v]{ac}e^{i(2h+1)\pi/(u+v)},\quad
h=0,1,\ldots,u+v-1.
\end{equation}
It can be easily shown that all of these roots lie in the sector
of stability defined by (\ref{sector_cond}) if and only if we have
\begin{equation}
\frac{\pi}{u+v}>\frac{\pi}{2M},
\end{equation}
or equivalently,
\begin{equation}\label{temp1}
\frac{u}{M}+\frac{v}{M}=\alpha+\beta<2.
\end{equation}
But, inequality (\ref{temp1}) is always satisfied since it is
assumed that both $\alpha$ and $\beta$ are between 0 and 1. To sum
up, $\mathbf{Y}_2^*$ is the stable equilibrium point of
(\ref{lotka1}) if and only if we have $ac>0$ provided that
$\alpha$ and $\beta$ are real numbers between 0 and 1. Note that
the stability of $\mathbf{Y}_2^*$ is also independent of the value
assigned to $b$. Figure \ref{fig_reg} summarizes the stability
properties of $\mathbf{Y}_1^*$ and $\mathbf{Y}_2^*$ in $a-c$
plane.
\subsection{Case 2: $\alpha$ and $\beta$ between 1 and 2}
The stability theorem presented in Section \ref{subsec_stab}
cannot directly be used to study the stability of the equilibrium
points of (\ref{lotka1}) when $\alpha$ and/or $\beta$ are greater
than unity. In order to study the stability of the equilibrium
points of (\ref{lotka1}) when both $\alpha$ and $\beta$ are
rational numbers between 1 and 2, first we should write this
equation in an equivalent form such that the order of
differentiation in all equations be between 0 and 1. Assuming
$y_3=dy_1/dt$, $y_4=dy_2/dt$, $\alpha=1+\alpha_1$, and
$\beta=1+\beta_1$, (\ref{lotka1}) can be written as
\begin{equation}\label{lotka2}
\left\{%
\begin{array}{ll}
    D^{\alpha_1}y_3=y_1(a-by_2)=f_1(y_1,y_2,y_3,y_4),\\
    D^{\beta_1}y_4=y_2(-c+by_1)=f_2(y_1,y_2,y_3,y_4), \\
    Dy_1=y_3=f_3(y_1,y_2,y_3,y_4), \\
    Dy_2=y_4=f_4(y_1,y_2,y_3,y_4),\\
\end{array}%
\right.
\end{equation}
where $\alpha_1=u/M$ and $\beta_1=v/M$ are rational numbers
between 0 and 1. System (\ref{lotka2}) has two equilibrium points
denoted as $\mathbf{Y}_1^*=(0,0,0,0)$ and
$\mathbf{Y}_2^*=(c/b,a/b,0,0)$ (the state vector of the system is
considered as $(y_3,y_4,y_1,y_2)$). The Jacobian matrices at
$\mathbf{Y}_1^*$ and $\mathbf{Y}_2^*$ are obtained as
\begin{equation}
\mathbf{J}_1=\left(%
\begin{array}{cccc}
  a & 0 & 0 & 0 \\
  0 & -c & 0 & 0 \\
  0 & 0 & 1 & 0 \\
  0 & 0 & 0 & 1 \\
\end{array}%
\right),
\end{equation}
and
\begin{equation}
\mathbf{J}_2=\left(%
\begin{array}{cccc}
  0 & -c & 0 & 0 \\
  a & 0 & 0 & 0 \\
  0 & 0 & 1 & 0 \\
  0 & 0 & 0 & 1 \\
\end{array}%
\right),
\end{equation}
respectively. At the first equilibrium point, (\ref{det_eq}) reads
\begin{multline}\label{charac3}
\det\left(\left(%
\begin{array}{cccc}
  \lambda^u & 0 & 0 & 0 \\
  0 & \lambda^v & 0 & 0 \\
  0 & 0 & \lambda^M & 0 \\
  0 & 0 & 0 &  \lambda^M \\
\end{array}%
\right)-\left(%
\begin{array}{cccc}
  a & 0 & 0 & 0 \\
  0 & -c & 0 & 0 \\
  0 & 0 & 1 & 0 \\
  0 & 0 & 0 & 1 \\
\end{array}%
\right)\right)\\=\left(\lambda^u-a\right)\left(\lambda^v+c\right)\left(\lambda^M-1\right)^2=0.
\end{multline}
Equation (\ref{charac3}) always has two roots at $\lambda=1$,
which do not lie in the sector of stability defined by
(\ref{sector_cond}). Therefore, $\mathbf{Y}_1^*$ is an unstable
equilibrium point for (\ref{lotka1}) for all values of $\alpha$
and $\beta$ between 1 and 2 regardless of the values assigned to
$a$, $b$, and $c$.

At the second equilibrium point, (\ref{det_eq}) concludes that
\begin{multline}\label{charac4}
\det\left(\left(%
\begin{array}{cccc}
  \lambda^u & 0 & 0 & 0 \\
  0 & \lambda^v & 0 & 0 \\
  0 & 0 & \lambda^M & 0 \\
  0 & 0 & 0 &  \lambda^M \\
\end{array}%
\right)-\left(%
\begin{array}{cccc}
  0 & -c & 0 & 0 \\
  a & 0 & 0 & 0 \\
  0 & 0 & 1 & 0 \\
  0 & 0 & 0 & 1 \\
\end{array}%
\right)\right)\\=\left(\lambda^{u+v}+ac\right)\left(\lambda^M-1\right)^2=0,
\end{multline}
which also has two roots at $\lambda=1$, outside the sector of
stability. Hence, $\mathbf{Y}_2^*$ is also an unstable equilibrium
point for (\ref{lotka1}) for all values of $\alpha$ and $\beta$
between 1 and 2 regardless of the special values assigned to $a$,
$b$, and $c$.

\section{Domain of attraction of stable equilibrium points}\label{sec_domain}
In the previous section we proved that the equilibrium points of
(\ref{lotka1}) are unstable for all values of $\alpha$ and $\beta$
between 1 and 2. In the following, we discuss on the domain of
attraction of the stable equilibrium points of (\ref{lotka1}) when
both $\alpha$ and $\beta$ are between 0 and 1. According to the
lack of analytical tools, most of the discussions in this section
are based on numerical calculations.

\subsection{Domain of attraction of $\mathbf{Y}_1^*=(0,0)$}
\begin{figure}
\begin{center}
\includegraphics[width=9.5cm]{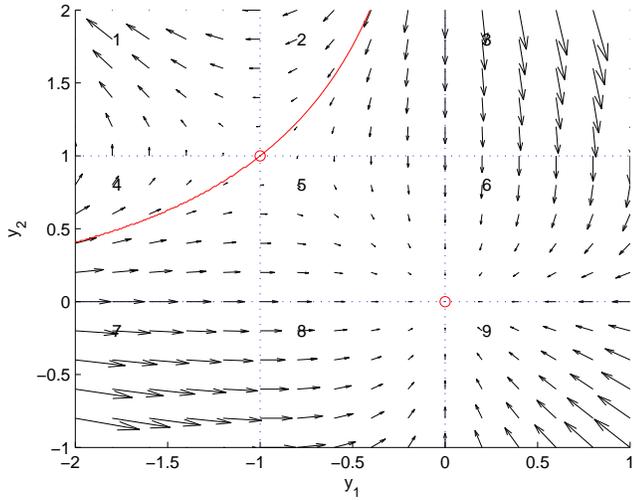}
\caption{Isoclines of (\ref{lotka1}) assuming $\alpha=\beta=-a =-
b=c=1$. The solid curve is the border of the domain of attraction
of $\mathbf{Y}_1^*$.}\label{fig_isoc1}
\end{center}
\end{figure}

According to Fig. \ref{fig_reg}, $\mathbf{Y}_1^*$ is a stable
equilibrium point for (\ref{lotka1}) (and consequently, has a
domain of attraction) if and only if we have $a<0$ and $c>0$,
provided that both $\alpha$ and $\beta$ are rational numbers
between 0 and 1. Figure \ref{fig_isoc1} shows the equilibrium
points of (\ref{lotka1}) and the corresponding isoclines in 9
different regions of the $y_1$-$y_2$ plane assuming $a = b= -c=-1$
(clearly, the following discussions can be extended to other
values of $a,~b,$ and $c$ provided that $a<0$ and $c>0$). In the
following, first we study the domain of attraction of
$\mathbf{Y}_1^*$ assuming $\alpha=\beta=1$ and then extend the
results to the case where both $\alpha$ and $\beta$ are between 0
and 1. Note that in the latter case the isoclines cannot directly
be used to determine the domain of attraction.


First, assume that $\alpha=\beta=1$. In this case, according to
the isoclines shown in Fig. \ref{fig_isoc1}, if $y_1(0)$ and
$y_2(0)$ lie in one of the
 regions 3, 5, 6, 7, 8, and 9 of Fig. \ref{fig_isoc1} then the
state trajectories of (\ref{lotka1}) move toward the origin by
increasing the time. But, if the initial conditions lie in region
1 then the state trajectories move toward infinity. It is also
clear that some parts of regions 2 and 4 in Fig. \ref{fig_isoc1}
belong to the domain of attraction of $\mathbf{Y}_1^*$. The border
of the domain of attraction in regions 2 and 4 can be determined
by solving the following initial value problem:
\begin{equation}
\frac{dy_2}{dy_1}=\frac{f_2}{f_1}=\frac{y_2(by_1-c)}{y_1(a-by_2)},\quad
y_2\left(\frac{c}{b}\right)=\frac{a}{b},
\end{equation}
the solution of which is implicitly given by
\begin{equation}
y_2^ay_1^c=\left(\frac{a}{b}\right)^a\left(\frac{c}{b}\right)^ce^{b(y_1+y_2)-(a+c)}.
\end{equation}
The solid curve in Fig. \ref{fig_isoc1} shows the solution of the
above equation. To sum up, all points in the right-hand side of
this curve belong to the domain of attraction of $\mathbf{Y}_1^*$
when $\alpha=\beta=1$. Note that regions 1, 5, 6, and 9 in Fig.
\ref{fig_isoc1} have the property that all of the state
trajectories that begin from any point inside them will remain in
these regions forever.

\begin{figure}
\begin{center}
\includegraphics[width=9.5cm]{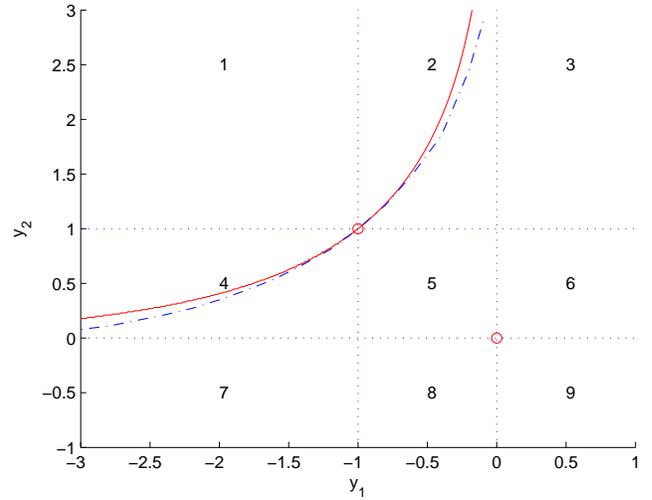}
\caption{The solid curve: the border of the domain of attraction
of $\mathbf{Y}_1^*$ when $\alpha=\beta=1$; the dash-dotted curve:
the border of the domain of attraction of $\mathbf{Y}_1^*$ when
$\alpha=\beta=0.1$.}\label{fig_dom1}
\end{center}
\end{figure}

In fractional case, however, according to the long memory
principle the border of the domain of attraction cannot be
determined analytically. The dash-dotted and solid curves in Fig.
\ref{fig_dom1} show the border of the domain of attraction when
$\alpha=\beta=0.1$ and $\alpha=\beta=1$, respectively assuming
$a=b=-c=-1$. In this figure, the dash-dotted curve has been
obtained by using numerical techniques (all simulations of this
paper are performed based on the numerical method proposed in
\citep{diethelm02}, which can be used to find the solution of a
Caputo definition based fractional differential equation).
Clearly, all points in the right-hand side of this curve belong to
the domain of attraction. As it can be observed, the domain of
attraction of $\mathbf{Y}_1^*$ becomes smaller by decreasing
$\alpha$ and $\beta$. Numerical simulations also show that the
border of the domain of attraction for the values of $\alpha$ and
$\beta$ between 0 and 0.1 is almost the same as the dash-dotted
curve in Fig. \ref{fig_dom1}. In general, precise simulations
confirm the fact that regions 3, 5, 6, 7, 8, and 9 in Fig.
\ref{fig_dom1} belong to the domain of attraction of
$\mathbf{Y}_1^*$ for all values of $\alpha$ and $\beta$ between 0
and 1. Moreover, all state trajectories that begin from any point
in regions 5, 6, 8, and 9 have the property that remain inside
these regions forever.

\subsection{Domain of attraction of $\mathbf{Y}_2^*=(c/b,a/b)$}
\begin{figure}
\begin{center}
\includegraphics[width=9.5cm]{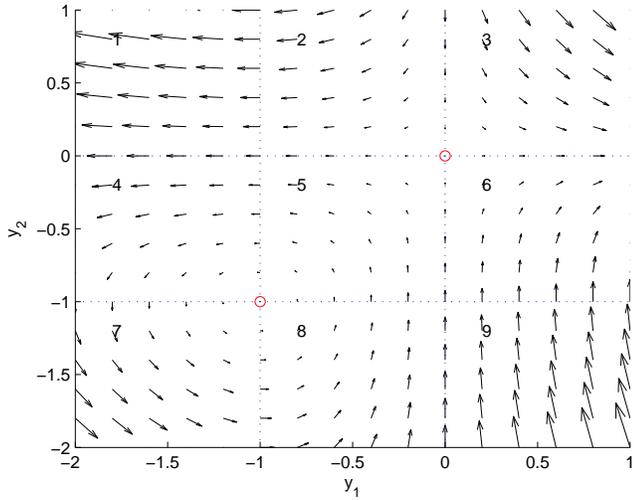}
\caption{Isoclines of (\ref{lotka1}) assuming $\alpha=\beta=a =
-b=c=1$.}\label{fig_isoc2}
\end{center}
\end{figure}

According to Fig. \ref{fig_reg}, $\mathbf{Y}_2^*$ is a stable
equilibrium point for (\ref{lotka1}) (and consequently, has a
domain of attraction) if and only if we have $a,c<0$ or $a,c>0$,
provided that both $\alpha$ and $\beta$ are rational numbers
between 0 and 1. Figure \ref{fig_isoc2} shows the equilibrium
points of (\ref{lotka1}) and the corresponding isoclines in 9
different regions of the $y_1$-$y_2$ plane assuming $a= b=-c=-1$.
According to this figure, the third quadrant of the $y_1$-$y_2$
plane is exactly equal to the domain of attraction of
$\mathbf{Y}_2^*$ when $\alpha=\beta=1$.

Similarly, precise numerical simulations show that the first
(third) quadrant of the $y_1$-$y_2$ plane is exactly equal to the
domain of attraction of $\mathbf{Y}_2^*$ when $b>0$ ($b<0$)
assuming that $a,c>0$ and $\alpha,\beta\in (0,1]$. It is also
observed that for all values of $a,c<0$ and $\alpha,\beta\in
(0,1]$ the first (third) quadrant of the $y_1$-$y_2$ plane is
exactly equal to the domain of attraction of $\mathbf{Y}_2^*$ when
$b<0$ ($b>0$).

Note that when $a$ and $c$ are chosen such that $\mathbf{Y}_2^*$
is a stable equilibrium point, it seems that it is a stable focus
for all values of $\alpha,\beta\in (0,1]$, but it will be shown in
the next section that this statement is not true.

\section{Some properties of the fractional-order Lotka-Volterra equation}
In this section, we numerically investigate some properties of the
two-dimensional fractional-order Lotka-Volterra equation assuming
$0<\alpha<1$ and $0<\beta<1$.

\begin{figure}
\begin{center}
\includegraphics[width=9.5cm]{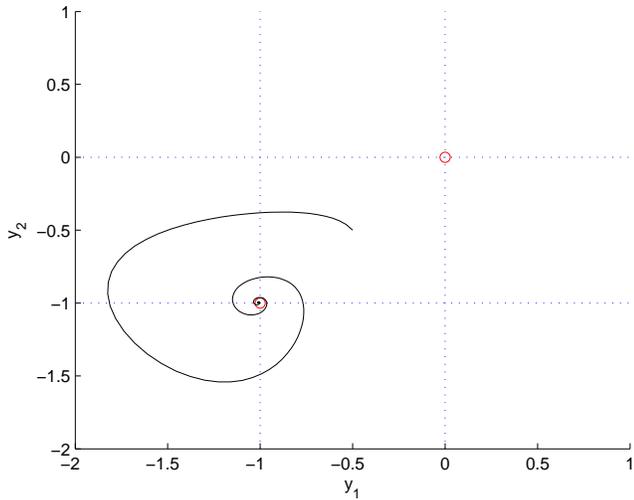}
\caption{The phase-plane portrait of (\ref{lotka1}) when
$\alpha=0.9$, $\beta=0.8$, $a = c = -b= 1$ and $y_1(0)= y_2(0) =
-0.5$ for $0\le t\le 80$.}\label{fig_simul1}
\end{center}
\end{figure}

\begin{figure}
\begin{center}
\includegraphics[width=9.5cm]{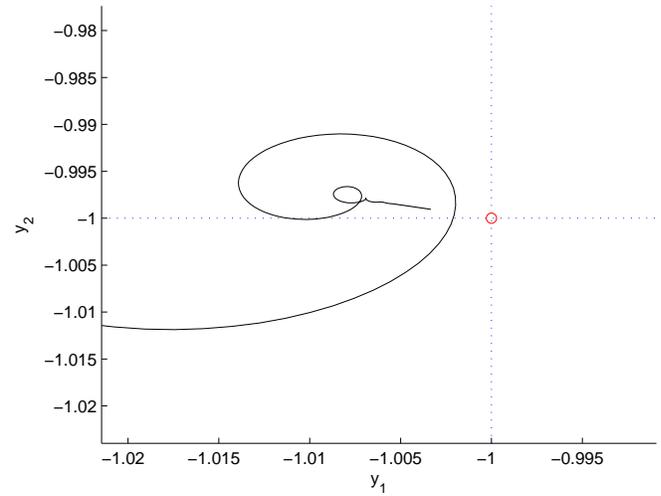}
\caption{Figure \ref{fig_simul1} focused around
$(-1,-1)$.}\label{fig_simul2}
\end{center}
\end{figure}

Figure \ref{fig_simul1} shows the phase-plane portrait of
(\ref{lotka1}) for $0\le t\le 80$ assuming $\alpha=0.9$,
$\beta=0.8$, $a = c = -b= 1$, and $y_1(0)= y_2(0) = -0.5$ . For
these values of parameters, the system has a stable equilibrium
point at $(-1,-1)$ and the third quadrant of the $y_1$-$y_2$ plane
is exactly equal to the domain of attraction of this equilibrium
point. As it is expected, the state trajectory of the system moves
toward this stable equilibrium point, which seems to be a stable
focus. Figure \ref{fig_simul2} shows the region around this stable
equilibrium point with more details. As it is observed, the state
trajectory of the system does not behave as it does near a stable
focus. Considerable number of simulations confirm the fact that
$\textbf{Y}_2^*$ (as well as $\textbf{Y}_1^*$) can never act as a
stable focus. As another fact, the state trajectory of the
fractional-order Lotka-Volterra system may intersect itself.

Note that inequality $D^{\alpha}y(t_0)>0,~(0<\alpha<1)$ does not
conclude that $y'(t_0)>0$, i.e. we may have $D^{\alpha}y(t_0)>0$
while $y$ is not increasing at $t=t_0$. That is why in the region
defined by $-2<y_2<-1$ and $-2<y_1<-1$ in Fig. \ref{fig_simul1},
$f_2$ is positive while $y_2$ is not a uniformly increasing
function of time.

\begin{figure}
\begin{center}
\includegraphics[width=9.5cm]{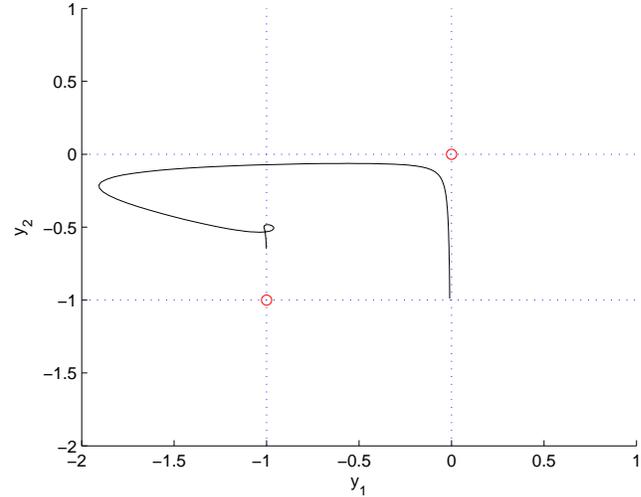}
\caption{The phase-plane portrait of (\ref{lotka1}) when
$\alpha=0.2$, $\beta=0.9$, $a = c = -b= 1$, $y_1(0)= -0.01$ and
$y_2(0) = -0.99$. }\label{fig_simul3}
\end{center}
\end{figure}

\begin{figure}
\begin{center}
\includegraphics[width=9.5cm]{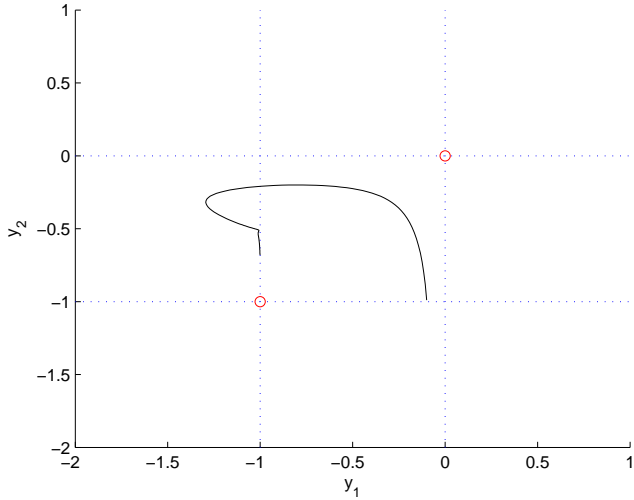}
\caption{The phase-plane portrait of (\ref{lotka1}) when
$\alpha=0.2$, $\beta=0.9$, $a = c = -b= 1$, $y_1(0)= -0.1$ and
$y_2(0) = -0.99$.}\label{fig_simul4}
\end{center}
\end{figure}

Another observation is that the qualitative behavior of the state
trajectories of (\ref{lotka1}), in general, depend on the initial
conditions of the system. For example, consider (\ref{lotka1})
with $\alpha=0.2$, $\beta=0.9$, and $a = c = -b= 1$ subject to the
initial conditions $y_1(0)= -0.01$ and $y_2(0) = -0.99$. The phase
portrait of this system is shown in Fig. \ref{fig_simul3} for
$0\le t\le 80$. As it is observed, the state trajectory of system
intersects itself and creates a \emph{tie}. Figure
\ref{fig_simul4} shows the phase portrait of the same system
assuming the initial conditions $y_1(0)=-0.1$ and $y_2(0)=-0.99$.
In this figure, the tie has been removed. This property is a
direct consequence of the long memory principle which exists only
in fractional-order systems.

Note that unlike the integer-order Lotka-Volterra equation which
exhibits stable limit cycles, the state trajectories of
(\ref{lotka1}) can never produce a limit cycle for the values of
$\alpha$ and $\beta$ between 0 and 1.

\section{Discussion and conclusion}
In this paper, we studied a general form of the two-dimensional
fractional-order Lotka-Volterra equation. Stability of the
equilibrium points was discussed and the domain of attraction of
stable equilibrium points was determined in several cases. It was
also numerically shown that the stable equilibrium points of this
system can never act as a stable focus. As another fact, it was
observed that the qualitative behavior of the state trajectories
of this system depend on the initial conditions. Specially, we
observed that for certain values of the parameters and initial
conditions the phase portrait of this system can intersect itself
and produce one or more ties.

Some other questions still remain unanswered. It is a well-known
fact that for a system described by
\begin{equation}\label{stat_eq}
\dot{\mathbf{x}}=\mathbf{f}(\mathbf{x},t),
\end{equation}
where $\mathbf{x}\in\mathbb{R}^n$ is the state vector of system,
the type of the equilibrium point $\mathbf{x}=\mathbf{x}^*$ can be
investigated by examining the eigenvalues of the Jacobian matrix
$\mathbf{J}=\partial
\mathbf{f}/\partial\mathbf{x}|_{\mathbf{x}=\mathbf{x}^*}$
\citep{slotine}. It concludes that $\mathbf{x}=\mathbf{x}^*$ is
either a stable (unstable) node or a stable (unstable) focus or a
saddle point. Note that in this case the eigenvalues of the
Jacobian matrix are actually the modes of the linearized system.
But, the discussions of this paper showed that according to the
long memory principle, the type of the equilibrium points of a
fractional-order system
 cannot be determined by
investigating the eigenvalues of the corresponding Jacobian
matrix. More precisely, linearizing (\ref{lotka1}) around any of
its equilibrium points leads to
\begin{equation}\label{lotka_lin}
\left(%
\begin{array}{c}
  D^{\alpha} Y_1 \\
  D^{\beta} Y_2 \\
\end{array}%
\right)=\left(%
\begin{array}{cc}
  a & -by_1^2 \\
  by_2^2 & -c \\
\end{array}%
\right)^*\left(%
\begin{array}{c}
  Y_1 \\
  Y_2 \\
\end{array}%
\right),
\end{equation}
where $Y_1=y_1-y_1^*$ and $Y_2=y_2-y_2^*$. By taking the Laplace
transform from both sides of the above equation, the
characteristic equation of the linearized system is obtained as
\begin{equation}
\left|\begin{array}{cc}
  s^\alpha-a & b(y_1^*)^2 \\
  -b(y_2^*)^2 & s^\beta+c \\
\end{array}\right|=0,
\end{equation}
where $s$ stands for the Laplace variable. The roots of the above
equation are modes of the linearized system (\ref{lotka_lin}),
which are distributed on a Riemann surface with infinity number of
Riemann sheets (clearly, if both $\alpha$ and $\beta$ be rational
numbers then the corresponding Riemann surface will have a finite
number of Riemann sheets). In fact, the discussions of this paper
showed that, unlike the classical case, in dealing with
fractional-order systems the modes of the linearized system cannot
be used the determine the type of the equilibrium points of
system. The relation between the modes of the linearized
fractional-order system and the type of its equilibrium points
remains as a challenging question.


\end{document}